\documentclass[fleqn]{mat01}
\usepackage{times,mathtimy,amssymb,latexsym}
\begin{document}

\setcounter{page}{175}
\firstpage{175}

\newtheorem{theor}{\bf Theorem}
\newtheorem{defini}[theor]{\rm DEFINITION}
\newtheorem{lem}[theor]{Lemma}
\newtheorem{rema}[theor]{Remark}
\newtheorem{exam}[theor]{Example}
\newtheorem{propo}[theor]{\rm PROPOSITION}
\newtheorem{coro}[theor]{\rm COROLLARY}
\newtheorem{case}{\it Case}

\font\xxxx=msam10 at 10pt
\def\Box{\mbox{\xxxx{\char'245}}}

\title{Quotient normed cones}

\markboth{Oscar Valero}{Quotient normed cones}

\author{OSCAR VALERO}

\address{Departamento de Ciencias Matem\'{a}ticas e
Inform\'{a}tica, Edificio Anselm Turmeda, Universidad de las Islas
Baleares, 07122 Palma, Spain\\
\noindent E-mail: o.valero@uib.es}

\volume{116}

\mon{May}

\parts{2}

\pubyear{2006}

\Date{MS received 16 December 2005}

\begin{abstract}
Given a normed cone $(X,p)$ and a subcone $Y,$ we construct and
study the quotient normed cone $(X/Y,\tilde{p})$ generated by $Y$.
In particular we characterize the bicompleteness of
$(X/Y,\tilde{p})$ in terms of the bicompleteness of $(X,p),$ and
prove that the dual quotient cone $((X/Y)^{*},\|\cdot
\|_{\tilde{p},u})$ can be identified as a distinguished subcone of
the dual cone $(X^{*},\|\cdot \|_{p,u})$. Furthermore, some parts
of the theory are presented in the general setting of the space
$CL(X,Y)$ of all continuous linear mappings from a normed cone
$(X,p)$ to a normed cone $(Y,q),$ extending several well-known
results related to open continuous linear mappings between normed
linear spaces.
\end{abstract}

\keyword{Normed cone; extended quasi-metric; continuous linear
mapping; bicom-\break pleteness; quotient cone; dual.}

\maketitle

\section{Introduction and preliminaries}

In recent years many works on functional analysis have been
obtained in order to extend the well-known results of the
classical theory of normed linear spaces to the framework of
asymmetric normed linear spaces and quasi-normed cones. In
particular, the dual of an asymmetric normed linear space has been
constructed and studied in \cite{LER}. In the same reference an
asymmetric version of the celebrated Alouglu theorem has been
proved (see also \cite{RSV2}). Several appropiate generalizations
of the structure of the dual of an asymmetric normed linear space
can be found in \cite{Scinica} and \cite{Dualcompl}. Hahn--Banach
type theorems in the frame of quasi-normed spaces have been given
in \cite{Ale,HBanach} and \cite{Tix}. In \cite{Ga} and
\cite{OsSa}, the completion of asymmetric normed linear spaces and
quasi-normed cones have been explored. An asymmetric version of
the Riesz theorem for finite dimension linear spaces can be found
in \cite{GTOP}.

It seems interesting to point out that quasi-normed cones and
other related `nonsymmetric' structures from topological algebra
and functional analysis, have been successfully applied, in the
last few years, to several problems in theoretical computer
science, approximation theory and physics, respectively (see
\S\S~11 and 12 of \cite{Ku2}, and also
\cite{Ali,GRS,Alicante,KeRo,RSa,RS2,RS,Tix,AppSci}).

The purpose of this paper is to show that it is possible to
generate in a natural way a quotient quasi-normed cone from a
subcone of a given quasi-normed cone. Actually, we analyse when
such quotient cones are bicomplete. We also construct and study
the dual cone of a quasi-normed cone and we prove that it can be
identified as the dual of a quotient cone. This is done with the
help of an appropiate notion of a `polar' cone.

\looseness 1 Throughout this paper, $\Bbb{R}^{+},$ $\Bbb{N}$ and
$\omega $ will denote the set of nonnegative real numbers, the set
of natural numbers and the set of nonnegative integer
numbers,\break respectively.

Recall that a {\it monoid} is a semigroup $(X,+)$ with neutral
element $0$.

According to \cite{KeRo}, a {\it cone} (on $\Bbb{R}^{+})$ is a
triple $(X,+,\cdot)$ such that $(X,+)$ is an abelian monoid, and
$\cdot $ is a mapping from $\Bbb{R}^{+}\times X$ to $X$ such that
for all $x,y\in X$ and $r,s\in \Bbb{R}^{+}$:
\begin{enumerate}
\renewcommand\labelenumi{(\roman{enumi})}
\leftskip .35pc
\item $r\cdot (s\cdot x)=(rs)\cdot x$;

\item $r\cdot (x+y)=(r\cdot x)+(r\cdot y)$;

\item $(r+s)\cdot x=(r\cdot x)+(s\cdot x)$;

\item $1\cdot x=x$;

\item $0\cdot x=0.$\vspace{-.4pc}
\end{enumerate}

A cone $(X,+,\cdot)$ is called {\it cancellative} if for all
$x,y,z\in X $, $z+x=z+y$ implies $x=y$.

Obviously, every linear space $(X,+,\cdot)$ can be considered as a
cancellative cone when we restrict the operation $\cdot $ to
$\Bbb{R}^{+}\times X.$

Let us recall that a {\it linear mapping} from a cone
$(X,+,\cdot)$ to a cone $(Y,+,\cdot)$ is a mapping $f\hbox{\rm :}\
X\rightarrow Y$ such that $f(\alpha \cdot x+\beta \cdot y)=\alpha
\cdot f(x)+\beta \cdot f(y)$.

A {\it subcone} of a cone $(X,+,\cdot)$ is a cone
$(Y,+|_{Y},\cdot|_{Y}) $ such that $Y$ is a subset of $X$ and
$+|_{Y}$ and $\cdot|_{Y}$ are the restriction of $+$ and $\cdot$
to $Y,$ respectively.

A {\it quasi-norm} on a cone $(X,+,\cdot)$ is a function
$q\hbox{\rm :}\ X\rightarrow \Bbb{R}^{+}$ such that for all
$x,y\in X$ and $r\in \Bbb{R}^{+}$:
\begin{enumerate}
\renewcommand\labelenumi{(\roman{enumi})}
\leftskip .35pc
\item $x=0$ if and only if there is $-x\in X$ and $q(x)=q(-x)=0;$

\item $q(r\cdot x)=rq(x);$

\item $q(x+y)\leq q(x)+q(y).$\vspace{-.4pc}
\end{enumerate}

If the quasi-norm $q$ satisfies: (i$'$) $q(x)=0$ if and only if
$x=0,$ then $q$ is called a {\it norm} on the cone $(X,+,\cdot).$

A quasi-norm defined on a linear space is called an {\it
asymmetric norm} in \cite{GRS}, \cite{Ga} and \cite{LER}.

Our main references for quasi-pseudo-metric spaces are \cite{FL}
and \cite{Ku2}.

Let us recall that a {\it quasi-pseudo-metric} on a set $X$ is a
nonnegative real-valued function $d$ on $X\times X$ such that for
all $x,y,z\in X$: (i) $d(x,x)=0;$ (ii) $d(x,z)\leq d(x,y)+d(y,z).$

In our context, by a {\it quasi-metric on} $X$ we mean a
quasi-pseudo-metric $d$ on $X$ that satisfies the following
condition: $d(x,y)=d(y,x)=0$ if and only if $x=y.$

We will also consider {\it extended quasi-{\rm (}pseudo-{\rm
)}metrics}. They satisfy the above three axioms, except that we
allow $d(x,y)=+\infty.$

Each extended quasi-pseudo-metric $d$ on a set $X$ induces a
topology $\mathcal{T}(d)$ on\thinspace $X$ which has as a base the
family of open $d$-balls $\{B_{d}(x,r)\hbox{\rm :}\ x\in X,$
$r>0\},$ where $B_{d}(x,r)=\{y\in X\hbox{\rm :}\ d(x,y)<r\}$ for
all $x\in X$ and $r>0.$

A ($n$ {\it extended{\rm )} quasi-{\rm (}pseudo-{\rm )}metric
space} is a pair $(X,d)$ such that $X$ is a set and $d$ is a ($n$
extended) quasi-(pseudo-)metric on $X.$

Similarly to \cite{Ko}, an extended quasi-metric $d$ on a cone
$(X,+,\cdot)$ is said to be {\it invariant} if for each $x,y,z\in
X$ and $r\in \Bbb{R}^{+},$ $d(x+z,y+z)=d(x,y)$ and
$d(rx,ry)=rd(x,y)$.

If $d$ is a ($n$ extended) quasi-metric on a set $X$, then the
function $d^{s}$ defined on $X\times X$ by $d^{s}(x,y)=\max
\{d(x,y),d(y,x)\}$ is a ($n$ extended) metric on $X.$ An extended
quasi-metric $d$ on a set $X$ is said to be {\it bicomplete} if
$d^{s}$ is a complete extended metric on $X.$

In \cite{OsSa} it was shown that it is possible to generate in a
natural way extended quasi-metrics from quasi-norms on
cancellative cones, extending to the well-known result that
establishes that each norm on a linear space $X$ induces a metric
on $X.$ In particular it was proved the following result:

\begin{propo}\label{ep}$\left.\right.$\vspace{.5pc}

\noindent Let $p$ be a quasi-norm on a cancellative cone
$(X,+,\cdot)$. Then the function $d_{p}$ defined on $X\times X$ by
\begin{equation*}
d_{p}(x,y) = \begin{cases} p(a), &\hbox{if} \ x\in X \ \hbox{and}
\ y\in x+X \ \hbox{with} \ y=x+a\\
+ \infty, &\hbox{otherwise}
\end{cases}
\end{equation*}
is an invariant extended quasi-metric on $X.$ Furthermore for each
$x\in X,$ $r\in \Bbb{R}^{+}\backslash \{0\}$ and each $\varepsilon
>0,$ $rB_{e_{p}}(x,\varepsilon)=rx+\{y\in X\hbox{\rm :}\ p(y)<r\varepsilon
\},$ and the translations are $\mathcal{T}(d_{p})$-open.
\end{propo}

The (extended) Sorgenfrey topology is obtained as a particular
case of the above construction.

The following well-known example will be useful later on. For each
$x,y\in \Bbb{R}$, let $u(x)=x\vee 0.$ Then $u$ is clearly a
quasi-norm on $\Bbb{R}$ whose induced extended quasi-metric is the
so-called (extended) upper quasi-metric $d_{u}$ on $\Bbb{R}.$

According to \cite{OsSa}, a quasi-normed cone is a pair $(X,p)$
where $X$ is a cancellative cone and $p$ is a quasi-norm on $X.$

\section{The quotient cone of a quasi-normed cone}

We can find several techniques to construct new normed linear
spaces from a given one in the literature (see, for example,
\cite{Rudin,Jameson} and \cite{Kelley2}). In particular one of
them consists of deriving quotient linear spaces from a given
linear space. Let us briefly recall this construction.

Let $X$ be a linear space and let $Y$ be a linear subspace of $X.$
For each $x\in X,$ let $[x]=x+Y.$ The family of such subsets is
called the quotient of $X$ by $Y$ and it is denoted by $X/Y.$ The
usual set-addition defines a sum on $X/Y$ by $x+Y+z+Y=x+z+Y.$
Moreover, the product by scalars is defined in $X/Y$ by $\lambda
\cdot (x+Y)=\lambda \cdot x+Y.$ Under these operations $X/Y$ forms
a linear space. Furthemore, if $\|\cdot\|$ is a norm defined on
$X$ and $Y$ is closed then $(E/Y,\|\cdot\|_{E/Y})$ admits a normed
linear structure, with $\|[x]\|_{E/Y}=\inf\{\|x+y\|\hbox{\rm
:}\break y\in Y\}$.

Next we introduce a general method for generating quotient spaces
from a quasi-normed cone $(X,p),$ which preserves the quasi-normed
cone structure, obtaining as a particular case of our construction
the above classical technique.

Let $X$ be a cancellative cone and let $Y$ be a subcone of $X$.
Denote by $G_{Y}$ the set $\{y\in Y\hbox{\rm :}\ -y\in Y\}.$
Clearly $G_{Y}$ is a nonempty set because of $0\in G_{Y}$. It is
easy to see that $G_{Y}$ admits a subcone structure, in fact it is
a linear space.

The following relation will allow us to construct a quasi-normed
cone with quotient space structure.

A pair of elements $x,z\in X$ are $\mathcal{R}$-related if there
exists $g\in G_{Y}$ such that $x+g=z.$ In this case, we write
$x\mathcal{R}y.$

This relation is obviously an equivalence relation, i.e.,
reflexive, symmetric and transitive. We will denote by $[x]$ the
equivalence class of $x$ with respect to $\mathcal{R}$. Thus, the
classes are given by $[x]=\{x+y\hbox{\rm :}\ y\in G_{Y}\}$ and the
family of them by $X/Y=\{[x]\hbox{\rm :}\ x\in X\}.$

\begin{propo}$\left.\right.$\vspace{.5pc}

\noindent The set $X/Y$ endowed with the usual sum and product of
equivalence classes $[x]+[z]:=[x+z]$ and $\lambda \cdot
[x]:=[\lambda \cdot x]$ where $x,z\in X$ and $\lambda \in
\Bbb{R}^{+},$ is a cancellative cone.
\end{propo}

\begin{proof}
Let us show that the sum of equivalence classes is well-defined.
Indeed, if $v\in [x]$ and $w\in [z]$ there are elements
$g_{1},g_{2}\in G_{Y}$ such that $x=v+g_{1}$ and $z=w+g_{2}.$ Then
$x+z=v+w+g_{1}+g_{2}.$ Whence $[v+w]=[x+z].$ The natural product
is also well-defined because $\lambda \cdot x=\lambda \cdot
v+\lambda \cdot g_{1},$ so that $[\lambda \cdot x]=[\lambda \cdot
v].$

Direct calculations show that $(X/Y,+,\cdot)$ is a cone with
neutral element $[0]=G_{Y}$. Finally, we show that $X/Y$ has the
cancellation law. Let $x,y,z\in X$ such that $[x]+[y]=[x]+[z].$
Then $[x+y]=[x+z].$ It follows that there exists $g\in G_{Y}$ with
$x+y=x+z+g.$ Since $X$ is cancellative, $y=z+g,$ so that $y\in
[z].$ Consequently $[y]=[z].$\hfill $\Box$
\end{proof}

\begin{defini}$\left.\right.$\vspace{.5pc}

\noindent{\rm Let $X$ be a cancellative cone and $Y$ be a subcone.
Then $X/Y$ is called {\it the quotient cone of} $X$ {\it by} $Y.$}
\end{defini}

\begin{rema}\label{inverse}
{\rm Note that if $[x]+[y]=[0]$ for any $x,y\in X,$ then
$[x+y]=[0]$ and, as a consequence there exists $g\in G_{Y}$ such
that $x+y=0+g.$ Hence $x+y+(-g)=0$. Therefore $-x=y+(-g)$ and
$-x\in [y].$}
\end{rema}

\begin{rema}{\rm
In case $X$ is a linear space and $Y$ is a linear subspace of $X$
we obtain as a particular case of our construction the classical
technique for generating quotient linear spaces from linear spaces
mentioned above, since $[0]=Y$ and $[x]=x+Y.$}
\end{rema}

\begin{exam}{\rm
Let $(\Bbb{R}^{2},+,\cdot)$ be endowed with the usual operations.
Let $Y$ be the subcone of $\Bbb{R}^{2}$ given by
$Y=\{(x,y)\hbox{\rm :}\ y\geq 0\}$. Obviously the subcone
$A=\{(x,0)\hbox{\rm :}\ x\in \Bbb{R}\}$ satisfies that $A\subset
Y$ and $G_{Y}=A.$ Thus the (cancellative) quotient cone
$\Bbb{R}^{2}/Y=\{[(x,y)]\hbox{\rm :}\ (x,y)\in \Bbb{R}^{2}\},$
where $[(x,y)]=(x,y)+A=\{(z+x,y)\hbox{\rm :}\ z\in \Bbb{R}\}.$}
\end{exam}

In the sequel, if $A$ is a subset of a quasi-metric space $(X,d),$
we will denote by $\bar{A}^{d}$ the closure of $A$ with respect to
$\mathcal{T(}d\mathcal{)}.$

The following property, whose proof is well-known and we omit,
will be useful later on.

\begin{propo}\label{qmclausura}$\left.\right.$\vspace{.5pc}

\noindent Let $(X,d)$ be a quasi-metric space and let $Y\subset
X$. Then{\rm ,} $x\in \bar{Y}^{d}$ if and only if $d(x,Y)=0,$
where $d(x,Y)=\inf\{d(x,y)\hbox{\rm :}\ y\in Y\}.$
\end{propo}

Next we denote by $d(Y,\cdot)$ the function defined as
$d(Y,x)=\inf \{d(y,x)\hbox{\rm :}\ y\in Y\}.$

In Proposition~\ref{pinf}, we construct a quasi-norm on $X/Y$ from
a quasi-normed cone $(X,p)$. For this, it is necessary to assume
that $G_{Y}$ is a closed subcone of $X$ with respect to
$\mathcal{T}(d_{p})$ (see Example~\ref{ejprenorma} below).

\begin{propo}\label{pinf}$\left.\right.$\vspace{.5pc}

\noindent Let $(X,p)$ be a quasi-normed cone and let $Y$ be a
subcone of $X$ such that $G_{Y}$ is closed in
$(X,\mathcal{T}(d_{p})).$ Then the pair $(X/Y,\hat{p})$ is a
quasi-normed cone{\rm ,} where the function $\hat{p}\hbox{\rm :}\
X/Y\rightarrow \Bbb{R}^{+}$ is defined by
\begin{equation*}
\hat{p}([x])=\inf \{p(x+y)\hbox{\rm :}\ y\in G_{Y}\}.
\end{equation*}
\end{propo}

\begin{proof}
It is clear that
\begin{equation*}
\hat{p}([0])=\inf \{p(y)\hbox{\rm :}\ y\in G_{Y}\}=0.
\end{equation*}
Let $x\in X$ such that $-x\in X$ and
$\hat{p}([x])=\hat{p}([-x])=0.$ Then,
\begin{equation*}
d_{p}(x,G_{Y})\vee d_{p}(G_{Y},x)\leq (d_{p})^{s}(x,0)=\max
\{p(x),p(-x)\}<+\infty.
\end{equation*}
Hence, given $\varepsilon >0,$ there exist $a,b\in X$ and
$g_{1},g_{2},h\in G_{Y}$ with

\begin{enumerate}
\renewcommand\labelenumi{(\roman{enumi})}
\leftskip .35pc
\item $g_{1}=x+b$ and $x=g_{2}+a,$

\item $p(b)<d_{p}(x,G_{Y})+\varepsilon,$

\item $p(-x+h)<\hat{p}([-x])+\varepsilon.$
\end{enumerate}

Since $d_{p}(x,g_{1})=p(b)$ we deduce that
\begin{equation*}
\hat{p}([-x])\leq p(-x+g_{1})=p(b)<d_{p}(x,G_{Y})+\varepsilon
\end{equation*}
and
\begin{equation*}
d_{p}(x,G_{Y})\leq p(-x+h)<\hat{p}([-x])+\varepsilon.
\end{equation*}
Then $\hat{p}([-x])\leq d_{p}(x,G_{Y})\leq \hat{p}([-x]).$
Therefore $d_{p}(x,G_{Y})=0$. Similarly it was showed that
$\hat{p}([x])\leq d_{p}(-x,G_{Y})\leq \hat{p}([x]).$

Since $G_{Y}$ is closed in $(X,\mathcal{T}(d_{p}))$ and by
Proposition~\ref{qmclausura}, $x,-x\in G_{Y}$. Thus $[x]=[0].$

Next we prove that $\hat{p}$ is homogeneous with respect to
nonnegative real numbers. We distinguish two cases:

\begin{case}{\rm $\lambda =0.$ Then the homogeneousness is immediately
obtained by definition of $\hat{p}$ and the fact that $X/Y$ is a
cone.}
\end{case}

\begin{case}{\rm $\lambda >0.$ Then
\begin{align*}
\hat{p}(\lambda \cdot [x]) &= \inf \{p(\lambda x+g)\hbox{\rm :}\
g\in G_{Y}\}=\inf \left\{\lambda p \left(x+\frac{g}{\lambda}
\right)\hbox{\rm :}\ g\in G_{Y}\right\} \\[.3pc]
&= \inf \{\lambda p(x+h)\hbox{\rm :}\ h\in G_{Y}\}=\lambda \inf
\{p(x+h)\hbox{\rm :}\ h\in G_{Y}\}\\[.3pc]
&= \lambda \hat{p}([x]).
\end{align*}

It remains to show the triangular inequality. Let $x,z\in X$ and
$\varepsilon >0.$ Then there exist $g,h\in G_{Y}$ such that
\begin{equation*}
p(x+g)<\hat{p}([x])+\varepsilon /2 \quad \hbox{and}\quad
p(z+h)<\hat{p}([z])+\varepsilon /2.
\end{equation*}
It follows that
\begin{align*}
\hat{p}([x]+[z]) &= \hat{p}([x+z])\leq p((x+z)+(g+h))\\[.3pc]
&\leq p(x+g)+p(z+h)<\hat{p}([x])+\hat{p}([z])+\varepsilon.
\end{align*}
We conclude that $\hat{p}([x]+[z])\leq \hat{p}([x])+\hat{p}
([z]).$}\hfill $\Box \ \ $
\end{case}
\end{proof}

The following example shows that the condition that the subcone
$G_{Y}$ is closed in $(X,\mathcal{T}(d_{p}))$ can not be omitted
in the statement of the preceding theorem.

\begin{exam}\label{ejprenorma}{\rm Consider
$X=(\Bbb{R}^{2},+,\cdot)$ endowed with the usual operations $+$
and $\cdot $. Let $A$ be the subcone of $\Bbb{R}^{2}$ given by
$A=\{(x,x)\hbox{\rm :}\ x\in \Bbb{R}\}$. Thus, it is clear that
$G_{A}=A.$ On the other hand, the function $p(x,y)=u(x)+u(y)$ is a
quasi-norm on $\Bbb{R}^{2},$ so that $(\Bbb{R}^{2},p)$ is a
quasi-normed cone. Furthermore, $(2,3)\in \bar{A}^{d_{p}}$ because
of $d_{p}\big((2,3),\big(2-\frac{1}{n},2-\frac{1}{n}
\big)\big)=p\big(-\frac{1}{n},-1-\frac{1}{n}\big)=u\big(-\frac{1}{n}\big)+u\big(-1-\frac{1}{n}\big)=0$
for all $n\in \Bbb{N}$ and, as a consequence, $G_{A}$ is not
closed in $(\Bbb{R}^{2},\mathcal{T}(d_{p})).$ Finally, the
function $\hat{p}$ satisfies that
$\hat{p}([(2,-3)])=\hat{p}([(-2,3)])=0.$ Whence $\hat{p}$ is a
prenorm on $\Bbb{R}^{2}/A$ but is not a quasi-norm.}
\end{exam}

\section{Bicomplete quasi-normed quotient cones}

Let us recall that a ($n$ extended) quasi-metric space $(Y,q)$ is
said to be a bicompletion of the (extended) quasi-metric space
$(X,d)$ if $(Y,q)$ is a bicomplete (extended) quasi-metric space
such that $(X,d)$ is isometric to a dense subspace of the
(extended) metric space $(Y,q^{s}).$ It is well-known that each
(extended) quasi-metric space $(X,d)$ has an (up to isometry)
unique bicompletion $(\tilde{X},\tilde{d})$ (see \cite{Di,Sa}).

In \cite{OsSa} it was introduced that the notion of bicomplete
quasi-normed cone and the construction of the bicompleteness was
given. Following \cite{OsSa}, a quasi-normed cone $(X,p)$ is said
to be {\it bicomplete} if the induced extended quasi-metric
$d_{p}$ is bicomplete. In connection with bicompleteness of
quasi-normed structures, some results for quasi-normed monoids and
asymmetric normed linear spaces may be found in \cite{RSV}
and\break \cite{Ga}.

In this section we characterize the bicompleteness of the
quasi-normed quotient cone in terms of the bicompleteness of the
original quasi-normed cone.

\begin{lem}\label{dom}
Let $(X,p)$ be a quasi-normed cone. Then $d_{\hat{p}}([x],[y])\leq
d_{p}(x,y)$ for all $x,y\in X.$
\end{lem}

\begin{proof}
We distinguish two cases:

\setcounter{case}{0}
\begin{case}{\rm $d_{p}(x,y)=+\infty.$ Then
$d_{\hat{p}}([x],[y])=+\infty,$ because otherwise we have that
there exists $[z]\in X/Y$ such that $[y]=[x]+[z].$ It follows that
$y=x+z+g$ for any $g\in G_{Y}.$ Hence we obtain that
$d_{p}(x,y)\leq p(z+g),$ a contradiction with the hypothesis.}
\end{case}

\begin{case}{\rm Now we suppose that $d_{p}(x,y)<+\infty.$ Then, given
$\varepsilon >0,$ there exists $c_{\varepsilon} \in X$ such that
$y=x+ c_{\varepsilon}$ and $p(c_{\varepsilon} )<d_{p}(x,y)+
\varepsilon $. Therefore $[y]=[x]+[c_{\varepsilon} ]$ with
$d_{\hat{p}}([x],[y])\leq \hat{p} ([c_{\varepsilon}
])<p(c_{\varepsilon} )<d_{p}(x,y)+\varepsilon.$}\hfill $\Box\ \
\!$
\end{case}\vspace{-1pc}
\end{proof}

\begin{rema}\label{cancellative}{\rm
Note that if $X$ is a cancellative cone, then $x=y+z$ and $y=x+w$
imply $w=-z$ for all $x,y,z,w\in X.$}
\end{rema}

The next theorem provides necessary and sufficient conditions for
bicompleteness of the quasi-normed quotient cone.

\begin{theor}[\!]\label{bicomplete}
Let $(X,p)$\ be a quasi-normed cone and let $Y$ be a subcone of
$X$ such that $G_{Y}$ is closed in $(X,\mathcal{T}(d_{p})).$ Then
$(X,p)$ is bicomplete if and only if $(G_{Y},p|_{G_{Y}})$ and
$(X/Y, \hat{p})$ are bicomplete.
\end{theor}

\begin{proof}
First we assume that $(X,p)$\ is bicomplete. Let $([x_{n}])_{n\in
\Bbb{N}}$ be a Cauchy sequence in the extended metric space
$(X/Y,(d_{\hat{p}})^{s}).$ Then there exists a subsequence
$([x_{n_{i}}])_{i\in \Bbb{N}}$ such that, given $\varepsilon >0,$
there is $n_{0}\in \Bbb{N}$ with
$(d_{\hat{p}})^{s}([x_{n_{i}}],[x_{n_{i+1}}])<\varepsilon $
whenever $n_{i}\geq n_{0}.$ It follows that there exist
$[a_{n_{i,i+1}}],[b_{n_{i,i+1}}]\in X/Y,$ such that
\begin{equation*}
\lbrack x_{n_{i}}]=[x_{n_{i+1}}]+[a_{n_{i,i+1}}], \quad
[x_{n_{i+1}}]=[x_{n_{i}}]+[b_{n_{i+1,i}}]
\end{equation*}
and $\max\{\hat{p}([a_{n_{i,i+1}}]),\hat{p}([ba_{n_{i+1,i}}])\}<
\varepsilon $ for all $n_{i}\geq n_{0}.$ By Remark
\ref{cancellative}, $[a_{n_{i,i+1}}]=-[b_{n_{i,i+1}}]$ and by
Remark \ref{inverse}, $-[a_{n_{i,i+1}}]=[-a_{n_{i,i+1}}].$
Consequently, there is $g_{n_{i,i+1}}\in G_{Y}$ which satisfies
\begin{equation*}
x_{n_{i}}=x_{n_{i+1}}+a_{n_{i,i+1}}+g_{n_{i,i+1}}, \quad
x_{n_{i+1}}=x_{n_{i}}+(-(a_{n_{i+1,i}}+g_{n_{i,i+1}}))
\end{equation*}
with
$\max\{p(a_{n_{i,i+1}}+g_{n_{i,i+1}}),p(-(a_{n_{i+1,i}}+g_{n_{i,i+1}}))
\}<\varepsilon $ whenever $n_{i}\geq n_{0}.$ Hence
$(x_{n_{i}})_{i\in \Bbb{N}}$ is a Cauchy sequence in the extended
metric space $(X,(d_{p})^{s}).$ Since $(X,p)$ is a bicomplete
quasi-normed cone there exists $x\in X$ such that
$\lim_{i\rightarrow +\infty} (d_{p})^{s}(x,x_{n_{i}})=0.$ Thus, by
Lemma~\ref{dom}, $\lim_{i\rightarrow +\infty}
(d_{\hat{p}})^{s}([x],[x_{n_{i}}])\leq \lim_{i\rightarrow +\infty
}(d_{p})^{s}(x,x_{n_{i}})=0.$ Therefore $([x_{n}])_{n\in \Bbb{N}}$
is a convergent sequence so that $(X/Y,\hat{p})$ is a bicomplete
quasi-normed cone. Furthermore, if $(X,p)$ is bicomplete we deduce
that all Cauchy sequence in $(G_{Y},(d_{p|_{G_{Y}}})^{s})$ are
convergent. Since the convergence with respect to
$\mathcal{T}(d_{p|_{G_{Y}}})^{s})$ implies convergence with
respect to $\mathcal{T}(d_{p|_{G_{Y}}})$ and $G_{Y}$ is closed in
$(X,\mathcal{T}(d_{p}))$ we obtain that $(G_{Y},p|_{G_{Y}})$
is\break bicomplete.

Conversely, let $(x_{n})_{n\Bbb{\in N}}$ be a Cauchy sequence in
$(X,(d_{p})^{s}).$ Given $\varepsilon >0,$ there exists $n_{0}\in
\Bbb{N}$ such that $d_{p}(x_{n},x_{m})<\varepsilon $ whenever
$n,m\geq n_{0}.$ Whence there exists $a_{nm}\in X$ such that
\begin{equation*}
x_{n}=x_{m}+a_{nm} \quad \hbox{and}\quad x_{m}=x_{n}+(-a_{nm})
\end{equation*}
with $\max\{p(a_{nm}),p(-a_{nm})\}<\varepsilon $ whenever $n,m\geq
n_{0}.$ By Lemma~\ref{dom} we obtain that
$(d_{\hat{p}})^{s}([x_{n}],[x_{m}])\leq
(d_{p})^{s}(x_{n},x_{m})=\max\{p(a_{nm}),p(-a_{nm})\}<\varepsilon
$ whenever $n,m\geq n_{0}.$ It follows that $([x_{n}])_{n\in
\Bbb{N}}$ is a Cauchy sequence in $(X/Y,(d_{\hat{p}})^{s}).$ Since
$(X/Y,\hat{p})$ is bicomplete, $\lim_{n\rightarrow \infty}
(d_{\hat{p}})^{s}([x],[x_{n}])=0$ for any $[x]\in X/Y.$ Hence
there exist $n_{1}\in \Bbb{N}$ and $[a_{n}]\in X/Y$ such that
$[x]=[x_{n}]+[a_{n}],$\ $[x_{n}]=[x]+[-a_{n}]$ and
$\max\{\hat{p}([a_{n}]),\hat{p}([-a_{n}])\}<\varepsilon$ whenever
$n\geq n_{1}.$ Put $n_{2}:=\max \{n_{0},n_{1}\}.$ Then there
exists $w_{n}\in G_{Y}$ such that $x_{n}=x+(-a_{n})+w_{n},
x_{m}=x+(-a_{m})+w_{m}$ and
$\max\{p((-a_{n})+w_{n}),p((-a_{m})+w_{m})\}<\varepsilon $
whenever $n,m\geq n_{2}$. Thus we deduce that
$w_{n}-a_{n}=w_{m}-a_{m}+a_{nm}$,
$w_{m}-a_{m}=w_{n}-a_{n}+(-a_{nm})$ and
$\max\{p(a_{nm}),p(-a_{nm})\}<\varepsilon $ whenever $m,n\geq
n_{2}.$ Whence we deduce that $d_{p}(w_{n}-a_{n},G_{Y})\leq
\hat{p}([a_{n}-w_{n}]=\hat{p}([a_{n}])<\varepsilon $ and
$d_{p}(G_{Y},w_{n}-a_{n})\leq
d_{p}(w_{n}-g,w_{n}-a_{n})=p(-a_{n}+g)<\hat{p}([-a_{n}])+\varepsilon
<2\varepsilon $ for all $n\geq n_{2}.$\break Since
\begin{align*}
\inf_{g\in G_{Y}}(d_{p})^{s}(g,w_{n}-a_{n}) &\leq \inf_{g\in
G_{Y}}\{d_{p}(w_{n}-a_{n},G_{Y})+d_{p}(G_{Y},w_{n}-a_{n})\}\\[.3pc]
&\leq \inf_{g\in G_{Y}}d_{p}(w_{n}-a_{n},G_{Y})+\inf_{g\in
G_{Y}}d_{p}(G_{Y},w_{n}-a_{n})\\[.3pc]
&= d_{p}(w_{n}-a_{n},G_{Y})+d_{p}(G_{Y},w_{n}-a_{n})<3\varepsilon,
\end{align*}
there exist $g_{n}\in G_{Y}$ and $t_{n}\in X$ such that
$g_{n}=w_{n}-a_{n}+t_{n}$, $w_{n}-a_{n}=g_{n}-t_{n}$ and
$\max\{p(t_{n}),p(-t_{n})\}<3\varepsilon $ whenever $n\geq n_{2}.$
Thus $g_{n}=w_{n}-a_{n}+t_{n}= w_{m}-a_{m}+a_{nm}+t_{n}=
g_{m}-t_{m}+a_{nm}+t_{n}$ and $g_{m}=g_{n}-t_{n}-a_{nm}+t_{m}$
with $\max \{p(-t_{m}+a_{nm}+t_{n}),
p(-t_{n}-a_{nm}+t_{m})\}<7\varepsilon$ for all $m,n\geq n_{2},$ so
\begin{equation*}
(d_{p})^{s}(g_{n},g_{m})=\max \{p(-t_{m}+a_{nm}+t_{n}),
p(-t_{n}-a_{nm}+t_{m})\}<5\varepsilon
\end{equation*}
whenever $n,m\geq n_{2}.$ It follows that $(g_{n})_{n\in \Bbb{N}}$
is a Cauchy sequence in the extended metric space
$(G_{Y},(d_{p})^{s}|_{G_{Y}}).$ Hence $\lim_{n\rightarrow \infty}
d_{p}(g,g_{n})=0$ for any $g\in G_{Y}$ because
$(G_{Y},p|_{G_{Y}})$ is bicomplete. Therefore
\begin{equation*}
\lim_{n\rightarrow \infty} d_{p}(g,w_{n}-a_{n})=\lim_{n\rightarrow
\infty} d_{p}(g,g_{n}-t_{n})=0
\end{equation*}
and
\begin{equation*}
\lim_{n\rightarrow \infty} d_{p}(x+g,x_{n})=\lim_{n\rightarrow
\infty} d_{p}(x+g,x+w_{n}-a_{n})=0.
\end{equation*}
This concludes the proof.\hfill $\Box$
\end{proof}

As a consequence of Theorem~\ref{bicomplete} we have the following
well-known result for normed linear spaces (see \cite{Conway}).

\begin{coro}$\left.\right.$\vspace{.5pc}

\noindent Let $(X,\|\cdot \|)$ be a normed linear space and let
$Y$ be a closed subspace of $X.$ Then $(X,\|\cdot \|)$\ is a
Banach space if and only if $(Y,\|\cdot \|)$\ and $(X/Y,\|\cdot
\|_{X/Y})$\ are Banach spaces.
\end{coro}

\section{Continuous linear mappings on quasi-normed cones}

Throughout this section we give several results concerning
continuous linear mappings, which will be useful in the next
section. A~mapping from a quasi-normed cone $(X,p)$ to a
quasi-normed cone $(Y,q)$ will be called {\it continuous} if it is
continuous from $(X,d_{p})$ to $(Y,d_{q})$.

In the sequel the set of all continuous linear mappings between
quasi-normed cones $(X,p)$ and $(Y,q)$ will be denoted by
$CL(X,Y).$ Obviously $CL(X,Y)$ is a cancellative cone for the
usual point-wise operations. In particular, when the pair $(Y,q)$
is the quasi-normed cone $(\Bbb{R},u\Bbb{)}$ we will denote
$CL(X,\Bbb{R})$ by $X^{*},$ and $X^{*}$ will be called the {\it
dual cone} of $X.$ Note that $f\in X^{*}$ if and only if it is a
linear and upper semicontinuous real-valued function on $(X,p).$
Some results concerning metrization and quasi-metrization for
several weak $^{*}$ topologies defined on the dual of a
quasi-normed cone and a generalization of Alouglu's theorem to
quasi-normed cone context were obtained in \cite{RSV2}.

A version of Theorem~\ref{continuity} and Proposition
\ref{Qnormed} for asymmetric normed linear spaces can be found in
\cite{Valentin} and \cite{LER}.

We omit the easy proof of the following lemma.

\begin{lem}\label{homo}
Let $(X,p)$ and $(Y,q)$ be two quasi-normed cones and let
$f\hbox{\rm :}\ X\rightarrow Y$ be a linear mapping. Then
$\frac{s}{r} f(B_{d_{p}}(0,r))=f(B_{d_{p}}(0,s))$ for all $r,s>0.$
\end{lem}

The next result characterizes the continuity of linear mappings
defined between quasi-normed cones. In order to help the reader we
include the complete proof of the below theorem although the
equivalence (1)--(3) has been proved in \cite{RSV2}.

\begin{theor}[\!]\label{continuity}
Let $(X,p)$ and $(Y,q)$\ be two quasi-normed cones and let
$f\hbox{\rm :}\ X\rightarrow Y$ be a linear mapping. Then the
following statements are equivalent{\rm :}
\begin{enumerate}
\renewcommand\labelenumi{\rm (\arabic{enumi})}
\leftskip .1pc
\item $f$ is continuous.

\item $f$ is bounded in $B_{d_{p}}(0,r)$ for every $r>0.$

\item There exists $c>0$ such that $q(f(x))\leq cp(x)$ for all $x\in
X.$\vspace{-.7pc}
\end{enumerate}
\end{theor}

\begin{proof}
First, we prove that $(1)$ implies $(2).$ Given $\varepsilon =1,$
there exists $r>0$ such that $d_{q}(0,f(x))=q(f(x))<1$ whenever
$d_{p}(0,x)=p(x)<r. $ Thus, $q(f(x))<1$ for all $x\in
B_{d_{p}}(0,r)$. Consider $r_{0}>0$. Next we show that $f$ is
bounded in $B_{d_{p}}(0,r_{0}).$ To this end we only consider that
$r_{0}>r$, otherwise there is nothing to prove. By
Lemma~\ref{homo}, $\frac{r_{0}}{r}
f(B_{d_{p}}(0,r))=f(B_{d_{p}}(0,r_{0})).$ It follows that
$q(f(y))\leq \frac{r_{0}}{r}$ for all $y\in B_{d_{p}}(0,r_{0}),$
and so $f$ is bounded. Since $f$ is bounded in $B_{d_{p}}(0,r)$
with $r>0$ and $B_{d_{q}}(0,0)\subseteq B_{d_{q}}(0,r),$ $f$ is
bounded in $B_{d_{q}}(0,0).$

Next we show that $(2)$ implies $(3)$. Let $M,r\in
\Bbb{R}^{+}\backslash \{0\}$ satisfying $q(f(x))\leq M$ for all
$x\in X$ such that $p(x)<r.$ Let $y\in X$. Now we show that
$q(f(y))\leq cp(y).$ We distinguish two cases.

\setcounter{case}{0}
\begin{case}{\rm
$p(y)=0.$ Suppose that $q(f(y))\neq 0$. Then $0<q(f(y))<M.$ Hence
there exists $\delta >0$ such that $\delta \cdot q(f(y))>M.$ It
follows that $q(f(\delta \cdot y))>M$ with $p(\delta \cdot
y)=\delta p(y)=0,$ which contradicts $(2)$.}
\end{case}

\begin{case}{\rm
$p(x)>0.$ Then
$p\big(\frac{r}{2}\frac{y}{p(y)}\big)=\frac{r}{2}\frac{p(y)}{p(y)}=\frac{r}{2}<r.$
Furthermore, by $(2),$
$q\big(f\big(\frac{r}{2}\frac{y}{p(y)}\big)\big)\leq M.$
Immediatelly we deduce that $\frac{r}{2}\frac{1}{p(y)}q(f(y))\leq
M.$ Therefore $q(f(y))\leq \frac{2}{r}M$ $p(y).$ Put
$c=\frac{2}{r}M.$ Then $q(f(y))\leq cp(y)$ for all $y\in X.$

Finally, we prove that $(3)$ implies $(1).$ Indeed, let
$\varepsilon >0,$ $\delta =\frac{\varepsilon}{c}$ and $x,y\in X$
such that $d_{p}(x,y)<\delta,$ then there exists $b\in X$ with
$y=x+b$ and $p(b)<\delta.$ Thus $d_{q}(f(x),f(y))\leq q(f(b))\leq
cp(b)<c\delta =\varepsilon.$ This concludes the proof.}\hfill
$\Box\ \ \!$
\end{case}\vspace{-1pc}
\end{proof}

\begin{rema}\label{product}{\rm
Let $(X,+,\cdot)$ be a cancellative cone. Then, $x=0$ if and only
if $\lambda \cdot x=0$ for any $\lambda >0.$}
\end{rema}

\begin{proof}
Suppose that $\lambda \cdot x=0$ for any $\lambda >0.$ Since
$\lambda \cdot x+\lambda \cdot y=\lambda \cdot y$ for all $y\in
X$, $\frac{1}{\lambda} (\lambda \cdot x)+\frac{1}{\lambda}
(\lambda \cdot y)=\frac{1}{\lambda} (\lambda \cdot y).$ Whence
$x+y=y$ for all $y,$ and so $x=0.$

Conversely, consider that $x=0.$ Then, for every $\lambda >0,$
$0+\frac{1}{\lambda} \cdot z=\frac{1}{\lambda} \cdot z$ for all
$z\in X.$ Thus $\lambda \cdot 0+z=z$ for all $z\in X.$ Whence
$\lambda \cdot 0=0$ for all $\lambda >0.$\hfill $\Box$
\end{proof}

Observe that $\overline{B_{d_{p}}}(x,r)\neq
\overline{B_{d_{p}}(x,r)}^{d_{p}},$ where
$\overline{B_{d_{p}}}(x,r)=\{y\in X\hbox{\rm :}\ d_{p}(x,y)\leq
r\}$ and $r\geq 0.$

\begin{exam}{\rm
Let $(\Bbb{R},+,\cdot)$ be the real numbers endowed with the usual
operations. It is clear that $(\Bbb{R},u)$ is a quasi-normed cone
such that $\overline{B_{d_{u}}}(0,1)=(-\infty,1]$ and
$\Bbb{R}\backslash \overline{B_{d_{u}}}(0,1)=(1,+\infty).$ So
$\overline{B_{d_{u}}}(0,1)$ is not a closed set in
$(\Bbb{R},\mathcal{T}(d_{u})).$}
\end{exam}

\begin{lem}\label{null}
Let $(X,p)$\ and $(Y,q)$\ be two quasi-normed cones and let $f\in
CL(X,Y)$. If $f|_{\bar{B}_{d_{p}}(0,1)}\equiv 0$\ then $f\equiv 0$
on $X.$
\end{lem}

\begin{proof}
Let $x\in X$. We distinguish two cases.

\setcounter{case}{0}
\begin{case}{\rm
$p(x)=0.$ Then $d_{p}(0,x)=p(x)=0$ and $x\in
\bar{B}_{d_{p}}(0,1),$ so $f(x)=0.$}
\end{case}

\begin{case}{\rm $p(x)\neq 0.$ Then $\frac{x}{p(x)}\in
\bar{B}_{d_{p}}(0,1)$ and $f\big(\frac{x}{p(x)}\big)=0$.
Consequently $\frac{1}{p(x)}$ $f(x)=0$ and by
Remark~\ref{product}, $f(x)=0.$}\hfill $\Box\ \ \!$
\end{case}\vspace{-1pc}
\end{proof}

\begin{propo}\label{Qnormed}$\left.\right.$\vspace{.5pc}

\noindent Let $(X,p)$ and $(Y,q)$ be two quasi-normed cones. Then
$(CL(X,Y),\|\cdot \|_{p,q})$\ is a quasi-normed cone{\rm ,} where
$\|f\|_{p,q}:=\sup_{x\in \overline{B_{d_{p}}}(x,1)}q(f(x)).$
\end{propo}

\begin{proof}
Obviously, the function $\|\cdot \|_{p,q}$ is well-defined by
Theorem~\ref{continuity}. Next, we only have to show that $\|\cdot
\|_{p,q}$ is a quasi-norm on $X.$ Indeed, $\|0\|_{p,q}=0.$ Let
$f\in CL(X,Y)$ satisfy $-f\in CL(X,Y)$ such that
$\|f\|_{p,q}=\|-f\|_{p,q}=0.$ Then
\begin{equation*}
\sup_{x\in \overline{B_{d_{p}}}(x,1)}q(f(x)) = \sup_{x\in
\overline{B_{p}}(x,1)} q(-f(x))=0.
\end{equation*}
Hence, we deduce that $f(x)=-f(x)=0$ for all $x\in
\overline{B_{d_{p}}}(0,1)$ and by Lemma~\ref{null}, $f=0.$

On the other hand, for each $\lambda \in \Bbb{R}^{+}$ and $f\in
CL(X,Y)$ we obtain
\begin{align*}
\|\lambda \cdot f\|_{p,q} &= \sup_{x\in
\overline{B_{d_{p}}}(x,1)}q(\lambda
f(x)) = \sup_{x\in \overline{B_{d_{p}}}(x,1)}\lambda q(f(x))\\[.3pc]
&= \lambda \sup_{x\in \overline{B_{d_{p}}}(x,1)}q(f(x))=\lambda
\|f\|_{p,q}.
\end{align*}

Finally, we show the triangle inequality. Let $f,g\in CL(X,Y).$
Then, given $\varepsilon >0,$ there exists $x_{\varepsilon} \in X$
such that $p(x_{\varepsilon} )\leq 1$ and
\begin{align*}
\|f+g\|_{p,q} &= \sup_{x\in
\overline{B_{d_{p}}}(x,1)}q(f(x)+g(x))<q(f(x_{\varepsilon})) +
q(g(x_{\varepsilon}))+\varepsilon\\[.3pc]
&\leq \|f\|_{p,q}+\|g\|_{p,q}+\varepsilon.
\end{align*}
Therefore $\|f+g\|_{p,q}\leq \|f\|_{p,q}+\|g\|_{p,q}.$\hfill
$\Box$
\end{proof}

We omit the easy proof of the following useful results.

\begin{propo}\label{normineq}$\left.\right.$\vspace{.5pc}

\noindent Let $(X,p)$\ and $(Y,q)$\ be two quasi-normed cones and
let $f\in CL(X,Y)$. Then $q(f(x))\leq \|f\|_{p,q}$\ $p(x)$ for all
$x\in X.$
\end{propo}

\begin{propo}\label{Three}$\left.\right.$\vspace{.5pc}

\noindent Let $(X,p),$\ $(Y,q)$\ and $(Z,w)$ be three quasi-normed
cones. If $f\in CL(X,Y)$ and $g\in CL(Y,Z)$ then $f\circ g\in
CL(X,Z)$ and $\|f\circ g\|_{p,w}\leq \|f\|_{p,q}\|g\|_{q,w}.$
\end{propo}

\begin{propo}\label{asynls}\label{injective}$\left.\right.$\vspace{.5pc}

\noindent Let $(X_{1},\|\cdot \|_{1})$ and $(X_{2},\|\cdot
\|_{2})$ be two asymmetric normed linear spaces and let
$T\hbox{\rm :}\ X_{1}\rightarrow X_{2}$ be a linear mapping. Then
the following assertions are equivalent{\rm :}

\begin{enumerate}
\renewcommand\labelenumi{\rm (\arabic{enumi})}
\leftskip .1pc
\item $T$ is injective and $T^{-1}$ is continuous.

\item There exists a positive constant $k$ such that $k\|x\|_{1}\leq
\|T(x)\|_{2}.$\vspace{-.5pc}
\end{enumerate}
\end{propo}

The key to the proof of Proposition~\ref{injective} is based on
the fact that a linear mapping is injective if and only if $\ker
T=\{x\in X\hbox{\rm :}\ T(x)=0\}=\{0\}$.

As a particular case of Proposition~\ref{injective} we obtain the
well-known result for the case of normed linear spaces (see, for
example \cite{Jameson}). However, Example~\ref{ker} shows that it
is possible to construct a noninjective mapping on the
quasi-normed cones which satisfies the condition (2) of the
preceding proposition. This is due to the fact that there are
functions whose $\ker$ is exactly the neutral element that are not
injective, so that the above proposition does not hold for
quasi-normed cones.

\begin{exam}\label{isometry}\label{ker}
{\rm Motivated by the applications to the analysis of complexity
of programs and algorithms given in \cite{Sch95}, it is introduced
and studied in \cite{RS} that the so-called dual complexity space,
which consists of the pair
$(\mathcal{C}^{*},d_{\mathcal{C}^{*}})$, where
\begin{equation*}
\mathcal{C}^{*}= \left\{f\in (\Bbb{R}^{+})^{\omega} \hbox{\rm :}\
\sum_{n=0}^{\infty }2^{-n}f(n)<+\infty \right\},
\end{equation*}
and $d_{\mathcal{C}^{*}}$ is the quasi-metric on $\mathcal{C}^{*}$
given by $d_{\mathcal{C}^{*}}(f,g)=\sum_{n=0}^{\infty}
2^{-n}[(g(n)-f(n))\vee 0]$.

Several properties of $d_{\mathcal{C}^{*}}$ are discussed in
\cite{RS}. In particular, observe that the topology induced by
$d_{\mathcal{C}^{*}}$ is not $T_{1}$.

On the other hand, $(\mathcal{C}^{*},+,\,\cdot \,)$ is clearly a
cancellative cone, with neutral element $f_{0}\in \mathcal{C}^{*}$
given by $f_{0}(n)=0$ for all $n\in \omega $, where $+$ is the
usual pointwise addition and $\cdot $ is the operation defined by
$(\lambda \cdot f)(n)=\lambda f(n)$ for all $n\in \omega$.

Let $p\hbox{\rm :}\ \mathcal{C}^{*}\rightarrow \Bbb{R}^{+}$
defined by $p(f)=\sum_{n=0}^{\infty} 2^{-n}f(n)$. It is routine to
see that $p$ is a quasi-norm on $\mathcal{C}^{*}$. Then the
induced extended quasi-metric $e_{p}$ on $\mathcal{C}^{*}$ is
given by
\begin{equation*}
e_{p}(f,g) = \begin{cases}
\sum_{n=0}^{\infty} 2^{-n} (g(n)-f(n)), &{\rm if} \ f\leq g\\
+\infty, &{\rm otherwise}
\end{cases}.
\end{equation*}
Several properties of the extended quasi-metric $e_{p}$ have been
studied in \cite{Alicante} and \cite{RSV}.

Let $X=\{f\in \mathcal{C}^{*}\hbox{\rm :}\ f(0)>0\}\cup
\{f_{0}\}$. It is routine to see that $X$ is a subcone of
$\mathcal{C}^{*}$.

Define $q\hbox{\rm :}\ X\rightarrow \Bbb{R}^{+}$ by $q(f)=f(0)$.
Clearly $q$ is a quasi-norm on $X$.

Let $F\hbox{\rm :}\ X\rightarrow \mathcal{C}^{*}$ defined by
$F(f)(0)=f(0)$ and $F(f)(n)=0 $ for all $f\in X$ and $n\in
\Bbb{N}$. Obviously $F$ is linear from $(X,+,\,\cdot \,)$ to
$(\mathcal{C}^{*},+,\,\cdot \,)$.

Moreover $p(F(f)) =\sum_{n=0}^{\infty} 2^{-n}F(f(n)) =f(0)=q(f)$
for all $f\in X$. So $F$ satisfies the condition $(2)$ of
Proposition~\ref{injective} and $\ker F=\{f_{0}\}.$ However, if
$f,g\in X$ satisfy $f(0)=g(0)$ and $f(1)\neq g(1)$, we obtain
$F(f)=F(g)$, and thus $F$ is not injective.}
\end{exam}

In order to obtain a quasi-normed cone version of
Proposition~\ref{injective} we introduce the next result.

\begin{propo}\label{finjective}$\left.\right.$\vspace{.5pc}

\noindent Let $(X,p)$\ and $(Y,q)$\ be two quasi-normed cones and
let $f\hbox{\rm :}\ X\rightarrow Y$\ be a linear mapping. If $f$\
is injective then the following statements are equivalent{\rm :}

\begin{enumerate}
\renewcommand\labelenumi{\rm (\arabic{enumi})}
\item There exists $c>0$ such that $p(f^{-1}(y))\leq cq(y)$ for
all $y\in f(X)$.

\item There exists $k>0$ such that $kp(x)\leq q(f(x))$ for all
$x\in X.$\vspace{-.5pc}
\end{enumerate}
\end{propo}
\pagebreak

\begin{proof}
First, we assume that the statement $(1)$ holds. Then for every
$x\in X$ we have $p(f^{-1}(f(x)))\leq c\,q(f(x))$ and
$kp(f^{-1}(f(x)))\leq q(f(x)),$ where $k=\frac{1}{c}$.

Now, assume that the statement $(2)$ holds and let $y\in f(X).$
Then there exists $x\in X$ such that $f(x)=y$ and $x=f^{-1}(y).$
Thus $kp(f^{-1}(y))\leq q(y),$ so that $p(f^{-1}(y))\leq cq(y),$
with $c=\frac{1}{k}.$\hfill $\Box$
\end{proof}

We finish the section extending some results related to open
continuous linear mappings between normed linear spaces to our
context. To this end, we introduce the $p$-{\it injectivity}.

\begin{defini}$\left.\right.$\vspace{.5pc}

\noindent{\rm Let $(X,p)$ be a quasi-normed cone and let $Y$ be a
nonempty subset of $X$. We will say that $f\hbox{\rm :}\
X\rightarrow Y$ is $p$-injective if $f(x)=f(y)$ implies
$p(x)=p(y)$ whenever $x,y\in X.$}
\end{defini}

Note that from the above definition one gathers that $\ker
f\subseteq \ker p$ and that all injective mappings are
$p$-injective$.$ Furthermore, if $(X,p)$ is an asymmetric normed
linear space the linearity joint with the $p$-injectivity implies
the injectivity. However, in the next example it is shown that the
equivalence is not true for quasi-normed cones.

\begin{exam}{\rm
Consider again, the quasi-normed cone $(X,q)$ of
Example~\ref{isometry}. It is easy to see that the mapping
$F\hbox{\rm :}\ X\rightarrow \mathcal{C}^{*}$ defined as in the
mentioned example is linear and $q$-injective, but it is not
injective.}
\end{exam}

\begin{theor}[\!]\label{open}
Let $(X,p)$ and $(Y,q)$ be two quasi-normed cones and let
$f\hbox{\rm :}\ X\rightarrow Y$ be a $p$-injective linear mapping.
Then the following assertions are equivalent{\rm :}

\begin{enumerate}
\renewcommand\labelenumi{\rm (\roman{enumi})}
\leftskip .35pc
\item $f$ is open and onto.

\item There exists $\delta >0$ such that $f(\overline{B_{d_{p}}}
(0,1))\supseteq \delta \overline{B_{d_{q}}}(0,1).$

\item There exists $M>0$ such that{\rm ,} given $y\in Y,$ there is $x\in X$
with $f(x)=y$ and $p(x)\leq Mq(y).$\vspace{-.6pc}
\end{enumerate}
\end{theor}

\begin{proof}
(i) clearly implies (ii). Next, suppose that (ii) holds and let
$y\in Y$. We distinguish two cases:

\setcounter{case}{0}
\begin{case}{\rm $q(y)=0.$ By hypothesis there exists $x\in
\overline{B_{d_{p}}}(0,1)$ such that $f(x)=\delta y$. Suppose that
$p(x)>0,$ then there exists $r>0$ with $rp(x)>1$ and
\begin{equation*}
f(r\cdot x)=r\cdot f(x)=r\cdot (\delta \cdot y),
\end{equation*}
since $q(r\cdot (\delta \cdot y))=r\delta q(y)=0,$ $r\cdot (\delta
\cdot y)\in \delta \overline{B_{d_{q}}}(0,1)$. On the other hand,
by (ii) we obtain that there exists $z\in
\overline{B_{d_{p}}}(0,1)$ such that $f(z)=\delta \cdot (r\cdot
y)=f(r\cdot x)$ with $r\cdot y\in \overline{B_{d_{q}}}(0,1).$
However $1<p(rx)$ and $p(z)\leq 1$ which contradicts the fact that
$f$ is $p$-injective$.$ Therefore $p(x)=0.$}
\end{case}

\begin{case}{\rm $q(y)\neq 0.$ Then, let
$y^{\prime}=\frac{\delta \cdot y}{q(y)}.$ By hypothesis, there
exists $x^{{\prime}}\in \overline{B_{d_{p}}}(0,1)$ such that
$f(x^{{\prime}})=y^{{\prime} }.$ Put
$x=\big(\frac{q(y)}{\delta}\big) \cdot x^{{\prime}}.$ Thus
$p(x)\leq q(y)/\delta.$

Next we prove that (iii) implies (i). Let $x\in X$ and let $V$ be
a neighborhood of $x.$ Then, by Proposition~\ref{ep}, there exists
$\varepsilon >0$ such that $x+z\in V$ for all $z$ with $p(z)\leq
\varepsilon$. Moreover, if $y\in Y$ satisfies $q(y)\leq
\varepsilon /M$ we have that there exists $z\in X$ such that
$f(z)=y$ and $p(z)\leq \varepsilon.$ Hence $f(x)+y=f(x+z)\in
f(V),$ and $f(V)$ is a neighborhood of $f(x).$ This concludes the
proof.}\hfill $\Box\ \ \!$\vspace{-1.5pc}
\end{case}
\end{proof}
\pagebreak

\begin{coro}$\left.\right.$\vspace{.5pc}

\noindent Let $(X,p)$ and $(X,q)$ be two quasi-normed cones. If
$f\hbox{\rm :}\ X\rightarrow Y$ is an open, bijective, continuous
linear mapping then it is an isomorphism.
\end{coro}

The condition of $p$-injectivity of Theorem~\ref{open} cannot be
omitted as seen in the following example.

\begin{exam}{\rm
Consider the cancellative cones $(\Bbb{R}^{2},+,\cdot)$ and
$(\Bbb{R},\pm,\odot)$ endowed with the usual operations. Define
$f\hbox{\rm :}\ \Bbb{R}^{2}\rightarrow \Bbb{R}$ by $f(x,y)=x.$ It
is clear that $f$ is linear and onto. Consider the quasi-norms
$p\hbox{\rm :}\ \Bbb{R}^{2}\rightarrow \Bbb{R}^{+}$ and
$u\hbox{\rm :}\ \Bbb{R}\rightarrow \Bbb{R}^{+}$ defined by
$p(x,y)=|y|$ and $u(x)=x\vee 0.$ Obviously $f(2,3)=f(2,4)$ and
$p(2,3)\neq p(2,4),$ and as a consequence $f$ is not
$p$-injective. Next we show that $f$ holds in (ii) in
Theorem~\ref{open}, but not in (iii). Indeed, we have that
$(-\infty,1]=\overline{B_{d_{u}}} (0,1)\Bbb{\subset}
f(\overline{B_{d_{p}}}(0,1))=\Bbb{R}$ and
$|y|=p(0,y)>Mu(f(0,y))=0,$ for all $M>0$ and $y\neq 0.$}
\end{exam}

\section{Quasi-normed cones as image of quotient cones and duality}

Given a normed linear space $(X,\|\cdot \|)$ and a linear subspace
$Y\subset X,$ it is defined the {\it polar of} $Y$ by
$Y^{0}=\{f\in X^{*}\hbox{\rm :}\ f|_{Y}\equiv 0\},$ where $X^{*}$
is the dual of $X.$ It is well-known that the polar inherits the
linear structure of $X^{*}$ and it is isometrically isomorphic to
the dual of the quotient linear space $X/Y$ whenever $Y$ is closed
(see \cite{Jameson})$.$ Next, we introduce a suitable notion of
polar which allows us to extend the mentioned property to our
context.

Define the mapping $\varphi \hbox{\rm :}\ X\rightarrow X/Y$ by
$\varphi (x)=[x].$ Obviously $\varphi $ is linear and onto$.$ Such
a mapping will be called the {\it quotient mapping} {\it of}
$X${\it \ induced by} $Y.$ Next we show that the mapping $\varphi
$ is, in addition, continuous.

\begin{propo}$\left.\right.$\vspace{.5pc}

\noindent Let $(X,p)$ be a quasi-normed cone and let $Y$ be a
subcone of $X$ such that $G_{Y}$ is closed in
$(X,\mathcal{T}(d_{p})).$ Then the quotient mapping $\varphi
\hbox{\rm :}\ X\rightarrow X/Y$\ is continuous.
\end{propo}

\begin{proof}
Let $x\in X.$ Then $\hat{p}(\varphi (x))=\hat{p}([x])= \inf_{y\in
G_{Y}}p(x+y)\leq \inf_{y\in G_{Y}}p(x)+p(y)=p(x).$\hfill
$\Box$\vspace{.3pc}
\end{proof}

\begin{coro}\label{phi}$\left.\right.$\vspace{.5pc}

\noindent $\|\varphi \|_{p,\hat{p}}\leq 1.$\vspace{.3pc}
\end{coro}

\begin{defini}$\left.\right.$\vspace{.5pc}

\noindent {\rm Let $(X,p)$ be a quasi-normed cone. We will call
{\it polar} of the subcone $Y$ of $X$ to the set $Y^{0}=\{f\in
X^{*}\hbox{\rm :}\ f|_{G_{Y}}\equiv 0 \dot{\}}$.}
\end{defini}

Note that as a particular case of our definition we have the polar
of a normed linear space, because of $G_{Y}=Y$ whenever $Y$ is a
linear subspace of a normed linear space.

Let $(X,p)$ and $(Y,q)$ be two quasi-normed cones. A~linear
mapping $T\hbox{\rm :}\ X\rightarrow Y$ is an {\it isometry} if
$q(T(x))=p(x)$ for all $x\in X.$ If, in addition, $T$ is an
isomorphism, i.e., $T$ is a continuous linear bijective mapping
with $T^{-1}$ continuous, then $T$ is called an {\it isometric
isomorphism}.

Observe that, contrary to the asymmetric linear case, there are
isometries on \hbox{(quasi-)} normed cones which are not injective
(see Example~\ref{ker} above). Two quasi-normed cones are called
{\it isometrically isomorphic} ({\it isomorphic}) if there exists
an isometric isomorphism (isomorphism) between them.

\begin{propo}$\left.\right.$\vspace{.5pc}

\noindent Let $(X,p)$ be a quasi-normed cone and let $Y$ be a
subcone of $X$\ such that $\overline{G_{Y}}^{d_{p}}=G_{Y}.$ Then
$Y^{0}$ is a subcone of $X^{*}$ isometrically isomorphic to
$(X/Y)^{*}.$
\end{propo}

\begin{proof}
We only show that $Y^{0}$ is isometrically isomorphic to
$(X/Y)^{*}$. Define the mapping $T\hbox{\rm :}\
(X/Y)^{*}\rightarrow Y^{0}$ by $(Tf)(x)=f(\varphi (x))=f([x]).$ It
is clear that $T$ is a well-defined mapping because:

\begin{enumerate}
\renewcommand\labelenumi{(\arabic{enumi})}
\leftskip .1pc
\item $(Tf)(g)=f([0])=0$ whenever $g\in G_{Y},$ since $f\in
(X/Y)^{*},$

\item $Tf$ is linear on $Y$ for all $f\in (X/Y)^{*},$ and
so is $\varphi$,

\item $Tf$ is continuous for all $f\in (X/Y)^{*}$, because
applying Proposition~\ref{normineq} for each $x\in X$ we have
\begin{equation*}
\hskip -1.25pc u((Tf)(x))=u(f(\varphi (x)))\leq \|f\|_{\hat{p},u}
\hat{p}(\varphi (x))\leq \|f\|_{\hat{p},u}\|\varphi
\|_{\hat{p},p}p(x).
\end{equation*}
\end{enumerate}

From the definition of $T$ we immediately deduce that it is
linear.

On the other hand, we obtain from $(3)$ and Corollary~\ref{phi}
that
\begin{equation*}
\|Tf\|_{p,u}=\sup \{u((Tf)(x))\hbox{\rm :}\ p(x)\leq 1\}\leq
\|f\|_{\hat{p},u}\|\varphi \|_{\hat{p},p}\leq \|f\|_{\hat{p},u}.
\end{equation*}
So $T$ is continuous.

Since $\tilde{p}([x])\leq p(x)$ for all $x\in X,$
$\|f\|_{\hat{p},u}\leq \|Tf\|_{p,u}.$ Hence
$\|f\|_{\hat{p},u}=\|Tf\|_{p,u}$ and $T$ is an isometry.

Next we show that $T$ is onto. For each $h\in Y^{0},$ define
$f_{h}\hbox{\rm :}\ X/Y\rightarrow \Bbb{R}$ as $f_{h}([x])=h(x).$
Thus $f_{h}$ is a well-defined function. Indeed,

\begin{enumerate}
\renewcommand\labelenumi{(\arabic{enumi})}
\leftskip .1pc
\item given $z\in [x]$ such that $z=x+g$ for any $g\in G_{Y}$ we
have
\begin{equation*}
\hskip -1.25pc f_{h}([z])=h(z)=h(x+g)=h(x)+h(g)=h(x)=f_{h}([x]).
\end{equation*}

\item $f_{h}$ is linear, and so is $h$.

\item Since for each $x\in X$ and $g\in G_{Y}$,
\begin{equation*}
\hskip -1.25pc u(f_{h}([x]))=u(h(x+g))\leq \|h\|_{p,u}p(x),
\end{equation*}
we follow that
\begin{equation*}
\hskip -1.25pc u(f_{h}([x]))\leq \|h\|_{p,u}\inf \{p(z)\hbox{\rm
:}\ z\in [x]\}=\|h\|_{p,u}\hat{p} ([x]),
\end{equation*}
and $f_{h}\in (X/Y)^{*}.$
\end{enumerate}

It is easy to see that $T(f_{h})=h.$

On the other hand, $T$ is injective. Suppose that $(Th)=(Tf)$ with
$f,h\in (X/Y)^{*}.$ Then $(Th)(x)=(Tf)(x)$ for all $x\in X.$
Whence $h([x])=f([x])$ for all $x\in X,$ and $f=h.$

Finally, we prove that $T^{-1}\hbox{\rm :}\ Y^{0}\rightarrow
(X/Y)^{*}$ is continuous. Let $h\in Y^{0},$ then
$T^{-1}(h)=f_{h}.$ It follows that
\begin{equation*}
u((T^{-1}h)([x]))=u(f_{h}([x]))\leq \|h\|_{p,u}\hat{p}([x]).
\end{equation*}
So $\|T^{-1}h\|_{\hat{p},u}\leq \|h\|_{p,u},$ and from
Theorem~\ref{continuity}, $T^{-1}$ is continuous.\hfill $\Box$
\end{proof}

In the classical context it is showed that, under open continuous
linear mappings, every normed linear space is isomorphic to a
distinguished quotient linear space (see \cite{Jameson}). To end
this section we obtain an asymmetric version of such a result.

Let $(X,p)$ and $(Y,q)$ be two quasi-normed cones and let
$T\hbox{\rm :}\ X\rightarrow Y$ be an onto linear mapping. Define
$\tilde{T}\hbox{\rm :}\ X/\ker T\rightarrow Y$ by
$\tilde{T}([x])=T(x).$ It is obvious that $\tilde{T}$ is
well-defined, since if $z\in [x]$, then $z=x+g$ for any $g\in
G_{\ker T}$ and $\tilde{T}([z])=T(z)=T(x)+T(g)=T(x).$

Since $T$ is linear and onto, so is $\tilde{T}$. Furthemore, if
$\varphi_{\ker T}$ is the quotient mapping induced by $\ker T,$ it
is a routine to check that $T=\tilde{T}\circ \varphi_{\ker T}.$

\begin{propo}\label{Ginject}$\left.\right.$\vspace{.5pc}

\noindent Let $(X,p)$ and $(Y,q)$ be two quasi-normed cones.
Then{\rm ,} $\tilde{T}$ is injective if and only if $T(x)=T(y)$
implies $y\in [x].$
\end{propo}

\begin{proof}
First, suppose that $\tilde{T}$ is injective. Then, if $T(x)=T(y)$
for any $x,y\in X$ we follow that $\tilde{T}([x])=\tilde{T}([y]).$
Thus $[x]=[y]$ and there exists $g\in G_{\ker T}$ such that
$y=x+g.$ So $T(y)=T(x).$

Now assume that $T(x)=T(y)$ implies $y\in [x].$ Then, if
$\tilde{T}([x])= \tilde{T}([y])$ for any $x,y\in X$ we have that
$y\in [x].$ Consequently $[y]=[x]$ and $\tilde{T}$ is
injective.\hfill $\Box$
\end{proof}

\begin{defini}$\left.\right.$\vspace{.5pc}

\noindent{\rm Let $(X,p)$ and $(Y,q)$ be two quasi-normed cones.
A~mapping $T\hbox{\rm :}\ X\rightarrow Y$ is $G$-injective if
$T(x)=T(y)$ implies $y\in [x]$.}
\end{defini}

Note that injectivity implies $G$-injectivity.

The below examples show that contrary to the case of normed linear
spaces \cite{Jameson}, $G_{_{\ker T}}$ is not closed in general in
$(X,\mathcal{T} (d_{p})).$

\begin{exam}{\rm
Consider $(\Bbb{R}^{2},+,\cdot)$ with the usual operations $+$ and
$\cdot$. Let $A$ be the subcone of $\Bbb{R}^{2}$ given by
$A=\{(x,0)\hbox{\rm :}\ x\in \Bbb{R}\}$ . Thus, it is clear that
$G_{A}=A.$ On the other hand, the function $p(x,y)=u(x)+u(y)$ is a
quasi-norm on $\Bbb{R}^{2},$ so that $(\Bbb{R}^{2},p)$ is a
quasi-normed cone. Furthermore, $(2,3)\in \bar{A}^{d_{p}}$ because
of $d_{p}((2,3),(-3n+2,0))=p((-3n,-3))=u(-3n)+u(-3)=0$ for all
$n\in \Bbb{N}$ and, as a consequence, $G_{A}$ is not closed in
$(\Bbb{R}^{2},\mathcal{T} (d_{p})).$ Define the linear mapping
$f(x,y)=y.$ Consequently $\ker f=A,$ $G_{_{\ker f}}= \ker f$ and
$G_{_{\ker f}}$ is not closed in $(X,\mathcal{T}(d_{p}))$.}
\end{exam}

\begin{exam}{\rm
Let $Y=\{(x,y)\hbox{\rm :}\ x\in \Bbb{R},y\geq 0\}$ be endowed
with the usual operations of $\Bbb{R}^{2}.$ Thus $(Y,p)$ is a
quasi-normed cone, where $p(x,y)=y.$ Define the linear mapping
$f\hbox{\rm :}\ Y\rightarrow \Bbb{R}$ by $f(x,y)=y.$ Hence $\ker
f=\{(x,0)\hbox{\rm :}\ x\in \Bbb{R}\}$ and $G_{_{\ker f}}=\ker f.$
Suppose that the sequence $\{(x_{n},0)\}_{n\in \Bbb{N}}\subset
G_{_{\ker f}}$ is convergent to $(x,y)$ in
$(Y,\mathcal{T}(d_{p})),$ then $(x_{n},0)=(x,y)+(z_{n},s_{n})$
eventually. It follows that $y=s_{n}=0$ eventually. Therefore
$(x,y)\in G_{_{\ker f}},$ and $G_{_{\ker f}}$ is closed in
$(Y,\mathcal{T}(d_{p})).$}
\end{exam}

In the light of the preceding examples we assume in the following
desired result that $\overline{G_{N}}^{d_{p}}=G_{N}.$

\begin{propo}$\left.\right.$\vspace{.5pc}

\noindent Let $(X,p)$ and $(Y,q)$ be two quasi-normed cones. Let
$T\hbox{\rm :}\ X\rightarrow Y$ be a onto{\rm ,} continuous linear
mapping such that $G_{_{\ker T}}$ is closed in
$(X,\mathcal{T}(d_{p})).$ Then the following holds{\rm :}

\begin{enumerate}
\renewcommand\labelenumi{\rm (\arabic{enumi})}
\leftskip .1pc
\item $\tilde{T}$\ is continuous.

\item $\|\tilde{T}\|_{\hat{p},q}=\|T\|_{p,q}.$

\item If in addition $T$ is $p$-injective{\rm ,} open and $G$-injective{\rm
,} then $\tilde{T}$ is an isomorphism.\vspace{-.7pc}
\end{enumerate}
\end{propo}

\begin{proof}
Let $[x]\in X/\ker f.$ By Proposition~\ref{normineq}, $q(\tilde{T}
([x]))=q(T(z))\leq \|T\|_{p,q}p(z)$ for all $z\in [x].$ Then
$q(\tilde{T} ([x]))\leq \|T\|_{p,q}\hat{p}([x]).$ Thus $\tilde{T}$
is continuous and $\|\tilde{T}\|_{\hat{p},q}\leq \|T\|_{_{p,q}}.$
Since $T=\tilde{T}\circ \varphi_{\ker f}$ we deduce that
$\|T\|_{p,q}=\|\tilde{T}\circ \varphi _{\ker f}\|_{p,q}$ and so,
by Proposition~\ref{Three}, $\|T\|_{p,q}\leq
\|\tilde{T}\|_{\hat{p},q}\|\varphi _{\ker f}\|_{p,\hat{p}}.$ By
Corollary~\ref{phi}, $\|T\|_{p,q}\leq \|\tilde{T}\|_{\hat{p},q}.$
We conclude that $\|T\|_{p,q}=\|\tilde{T}\|_{\hat{p},q}.$ Thus we
have proved (1) and (2).

Item (3) remains. Assume that $T$ is, in addition, open and
$p$-injective. Then, by Theorem~\ref{open}, there exists $M>0$
such that every element $y\in Y$ is expressible as $T(x)$ for any
$x\in X$, holding $p(x)\leq Mq(y).$ Since $T$ is $G$-injective, by
Proposition~\ref{Ginject}, $\tilde{T}$ is injective and
$\tilde{T}^{-1}(y)=[x].$ Furthermore, $\hat{p}([x])\leq p(x)\leq
Mq(y)$ and, by Proposition~\ref{finjective}, $\tilde{T}^{-1}$ is
continuous with $\|\tilde{T}\|_{\hat{p},q}\leq M.$\hfill $\Box$
\end{proof}

\section*{Acknowledgements}

The author acknowledges the support of the Spanish Ministry of
Science and Technology, Plan Nacional I+D+I, and FEDER, grant
BMF2003-02302.

\end{document}